\documentclass[11pt]{article}
\usepackage{amsmath}
\usepackage{CJK}
\usepackage{amsfonts}
\usepackage[dvips]{color}
\usepackage{amssymb}
\usepackage{ leftidx}
\usepackage{mathdots}
\usepackage{graphicx}
\usepackage{amssymb}
\usepackage{booktabs}

\makeatletter

\newcommand{\Rmnum}[1]{\expandafter\@slowromancap\romannumeral #1@}
\makeatother
\usepackage{fancyhdr}
\usepackage{amsfonts,latexsym,lscape}
\newtheorem{theorem}{Theorem}[section]

\newtheorem{lemma}{Lemma}[section]
\newtheorem{corollary}{Corollary}[section]

\newfont{\Bb}{msbm10 scaled\magstep1}

\newcommand{\R}{{\mathbf R}}

\def\sqr#1#2{{\vcenter{\hrule height .#2pt
      \hbox{\vrule width .#2pt height#1pt \kern#1pt\vrule width.#2pt}
                       \hrule height.#2pt}}}

\def\L{{\cal L}}

\def \[{\begin{equation}}
\def \]{\end{equation}}

\def\L{{\cal L}}

\def\bdes{\begin{description}}
\def\edes{\end{description}}
\def\benu{\begin{enumerate}}
\def\eenu{\end{enumerate}}
\def\bitm{\begin{itemize}}
\def\eitm{\end{itemize}}

\def\L{{\cal L}}

\def\R{{\sl I\kern-3.2pt R}}

\def\r{\right}

\def\sqr#1#2{{\vcenter{\hrule height .#2pt
      \hbox{\vrule width .#2pt height#1pt \kern#1pt\vrule width.#2pt}
                       \hrule height.#2pt}}}

\def\r1{{\mathbb{R}}}

\title{Stochastic Approximation Proximal Method of Multipliers  for Convex Stochastic Programming\thanks{Supported by the National Natural Science Foundation of China under project
No.11971089,  No.11731013 and No.11571059.}}
\author{Liwei Zhang\footnote{Institute of Operations Research and Control Theory, School of Mathematical Sciences, Dalian University
of Technology, Dalian 116024, China.(lwzhang@dlut.edu.cn)},
\,\,
Yule Zhang\footnote{Institute of Operations Research and Control Theory,
School of Mathematical Sciences, Dalian University of Technology,
Dalian 116024, China.(358768797@qq.com)},
\,\, Jia Wu\footnote{Institute of Operations Research and Control Theory,
School of Mathematical Sciences, Dalian University of Technology,
Dalian 116024, China.(wujia@dlut.edu.cn)}
}

\date{\today}
\begin{document}
\maketitle
 \vspace{2mm}
\begin{center}
\parbox{13.5cm}{\small
\textbf{Abstract.}
This paper considers the problem of minimizing a convex expectation function over a closed convex set, coupled with a set of  inequality  convex expectation constraints. We present a new stochastic approximation  type algorithm, namely the stochastic approximation proximal method of multipliers  (PMMSopt) to solve this convex stochastic optimization problem.
 We analyze regrets of a stochastic approximation proximal method of multipliers for solving convex stochastic optimization problems. Under mild conditions, we show that this algorithm exhibits  ${\rm O}(T^{-1/2})$  regret for both objective reduction and constraint violation if parameters in the algorithm are properly chosen, when the objective and constraint functions are generally convex, where $T$ denotes the number of iterations.  Moreover, we show that,
    with at least $1-e^{-T^{1/4}}$ probability, the algorithm has no more than  ${\rm O}(T^{-1/4})$ objective regret and no more than  ${\rm O}(T^{-1/8})$ constraint violation  regret.
  To the best of our knowledge, this is the first  time that such a proximal point method for solving expectation constrained stochastic optimization is presented in the literature.
\\[10pt]
 }
\end{center}

\section{Introduction}
\setcounter{equation}{0}
\quad \, In this paper, we consider the following stochastic optimization problem
\begin{equation}\label{eq:1}
\begin{array}{rl}
\displaystyle \min_{x \in X_0} & f(x)=\mathbb{E}[F(x,\xi)]\\[4pt]
{\rm s.t.} & g_i(x)=\mathbb{E}[G_i(x,\xi)] \leq 0, i=1,\ldots,p.\\
\end{array}
\end{equation}
Here $X_0 \subset \Re^n$ is a nonempty bounded closed convex set, $\xi$ is a random vector whose probability distribution $P$  is supported on set $\Xi \subseteq \Re^q$ and $F: {\cal O}_0 \times \Xi \rightarrow \Re$, $G_i:{\cal O}_0 \times \Xi \rightarrow \Re$, $i=1,\ldots, p$, where ${\cal O}_0\subset \Re^n$ is an open bounded convex set containing $X_0$.  Let $\Phi$ be the feasible region of Problem (\ref{eq:1}):
$$
\Phi=\left\{x\in X_0: g_i(x) \leq 0, i=1,\ldots,p\right\}.
$$
We assume that expectations
$$
\mathbb{E}[F(x,\xi)]=\displaystyle \int_{\Xi} F(x,\xi)dP(\xi),\, \mathbb{E}[G_i(x,\xi)]=\displaystyle \int_{\Xi} G_i(x,\xi)dP(\xi), i=1,\ldots, p
$$
are well defined and finite valued for every $x\in {\cal O}_0$.
Moreover, we assume that the  functions $F(\cdot,\xi)$ and
$G_i(\cdot,\xi)$, $i=1,\ldots,p$ are continuous and convex on
${\cal O}_0$. Denote $G(x,\xi)=(G_1(x,\xi),\ldots, G_p(x,\xi))^T$
and $g(x)=(g_1(x),\ldots, g_p(x))^T$, then
$$
g(x)=\displaystyle \int_{\Xi} G(x,\xi)dP(\xi).
$$

It is well-known that a computational difficulty of solving stochastic optimization problem (\ref{eq:1}) is that expectation is an multidimensional integral and it cannot be computed with a high accuracy for large dimension $q$. The aim of this paper is to construct a stochastic approximation proximal method of multipliers for solving  Problem (\ref{eq:1}). To this end we make the following assumptions.
\begin{itemize}
\item[(A1)] It is possible to generate an i.i.d. sample $\xi_1,\xi_2,\ldots,$ of realizations of random vector $\xi$.
    \item[(A2)] There is an oracle, which, for any point $(x,\xi)\in {\cal O}_0 \times \Xi$ returns stochastic subgradients $v_0(x,\xi),v_1(x,\xi),\ldots, v_p(x,\xi)$ of $F(x,\xi)$,$G_1(x,\xi),\ldots,G_p(x,\xi)$ such that
        $$
        v_i(x)=\mathbb{E}[v_i(x,\xi)],i=0,1,\ldots,p
        $$
        are well defined and are subgradients of $f(\cdot), g_1(\cdot),\ldots, g_p(\cdot)$ at $x$, respectively, i.e., $v_0(x)\in \partial f(x)$, $v_i(x) \in \partial g_i(x)$,$i=1,\ldots,p$.
        \item[(A3)] Let $D_0>0$, $\nu_f>0$ and $\nu_g>0$ such that
        $$ \|x'-x''\|\leq D_0, \forall x',x'' \in X_0
        $$
        and
        $$
        F(x',\xi)-F(x'',\xi) \leq \nu_f, \|G(x,\xi)\| \leq \nu_g, \forall x',x'' \in {\cal O}_0, \xi \in \Xi.
        $$
        \item[(A4)]
        Let $\kappa_f>0$ and $\kappa_g>0$ such that
                $$
        \|v_0(x,\xi)\| \leq \kappa_f, \,\, \|v_i(x,\xi)\| \leq \kappa_g,i=1,\ldots, p, \forall (x,\xi)\in {\cal O}_0\times \Xi,
        $$
        where $v_0(x,\xi)$ is a stochastic subgradient of $F(x,\xi)$ and $v_i(x,\xi)$ is a stochastic subgradient of $G_i(x,\xi)$,$i=1,\ldots,p$, $(x,\xi) \in {\cal O}_0 \times \Xi$.
        \item[(A5)] There exist $\epsilon_0>0$ and $\widehat x \in X_0$ such that
        $$
        g_i(\widehat x) \leq -\epsilon_0, \,\, i =1,\ldots, p.
        $$
            \end{itemize}
          The stochastic approximation (SA) technique  is going back to the pioneering paper by Robbins and Monro \cite{RM21}. Since then SA algorithms, due to low demand for computer memory and cheap computation cost at every iteration, became widely used in stochastic optimization and online optimization, see, e.g. \cite{P20},  \cite{Xiao2010} and \cite{Shai2011}. In the classical analysis of the SA algorithm,  initiated from the works \cite{C5} and \cite{S23},  it is assumed that $f(\cdot)$ is twice continuously differentiable and strongly convex and in the case when the minimizer of $f$ belongs to the interior of $\Phi$, exhibits asymptotically optimal rate of convergence $\mathbb{E}[f(x^t)-f^*]={\rm O}(t^{-1})$, where $x^t$ is $t$-th iterate and $f^*$ is the minimal value of $f(x)$ over $x\in \Phi$. This algorithm, however, is very sensitive to a choice of the respective stepsizes.  For overcoming this drawback, an important improvement of the SA method was developed by Polyak \cite{P18} and Polyak and Juditsky  \cite{P19}, where longer stepsizes were suggested with consequent averaging of the obtained iterates. Adopting the averaging technique to iterates generated by (under our notations)
          \begin{equation}\label{eq:projIt}
          x^{j+1}=\Pi_{\Phi}(x^j-\gamma_j v_0(x^j,\xi_j)),
          \end{equation}
    Nemirovski, Juditsky,  Lan, and  Shapiro \cite{Lan2009} shows that, without assuming smoothness and strong convexity of the objective function, the convergence  rate is ${\rm O}(t^{-1/2})$.  This paper  also demonstrates that a properly modified SA approach can be competitive and even significantly outperform the SAA method for a certain class of convex stochastic problems. After the seminal work of \cite{Lan2009}, there are many significant results appeared, even for non-convex stochastic optimization problems, see \cite{Lan2012a},\cite{Lan2012b},\cite{GLan2012},\cite{GLan2013},
 \cite{GLan2013a} and  \cite{GLZh2016}.
 Among all mentioned works, the feasible region set is an abstract closed convex set, none of these SA algorithms are applicable to expectation constrained problems. The computation of  projection $\Pi_{\Phi}$ is quite easy only when $\Phi$ is of a simple structure.
 However, when $\Phi$ is defined by  (\ref{eq:1}), the computation of projection $\Pi_{\Phi}$ is prohibitive, this is a difficult work. Therefore, it is quite important to obtain a numerical method with lower iteration complexity of both objective and constraint violation.

For stochastic optimization problems with expectation constraints,
Yu et. al \cite{YMNeely2017} proposed an algorithm
that achieves ${\rm O}(T^{-1/2})$ expected regret and constraint violations and ${\rm O}(T^{-1/2}\log T)$
high probability regret and constraint violations. Lan and Zhou \cite{LanZ2016} proposed a cooperative stochastic
approximation type algorithm inspired by Polayk's subgradient
method, which exhibits ${\rm O}(T^{-1/2})$ convergences in terms of
objective value and constraint violation.
In \cite{Boob}, the authors presented a novel Constraint Extrapolation
method for solving convex functional constrained problems which
includes the stochastic optimization problem (\ref{eq:1}) as a
special case. The method is a single-loop primal-dual type method,
which utilizes linear approximations of the constraint functions to
define the extrapolation step and also exhibits ${\rm O}(T^{-1/2})$
convergence results.

 A natural way to handle constraints for constrained optimization problems is to use augmented Lagrangian, which results in proximal point methods. It is well-known that Rockafellar \cite{rockafellar76B} proposed three proximal point methods for  convex programming, namely the proximal point method developed in \cite{rockafellar76A} applied to  maximum monotone inclusions of the primal optimality, the dual optimality and the saddle point optimality. The augmented Lagrange method, just the proximal point method applied to the dual optimality, has been studied deeply not only for convex optimization  but also for non-convex optimization. And the proximal point method for
  the primal optimality has been extensively implemented for solving various structured convex optimization problems. However,  the  proximal point method for the saddle point optimality, the so-called proximal method of multipliers by Rockafellar \cite{rockafellar76B}, has not been paid much attention.

  In this paper, we study the proximal method of multipliers for stochastic convex  problem, and analyze its regret bounds as well as probability guarantee for both objective reduction and constraint violation.

  For i.i.d. sample $\xi_1,\xi_2,\ldots$ of realizations of random vector $\xi$. Consider the following convex optimization problem
\begin{equation}\label{eq:Pt}
\begin{array}{ll}
\min & F(x,\xi_t)\\[4pt]
{\rm s.t.} & G_i(x,\xi_t) \leq 0,i=1,\ldots,p.\\
\end{array}
\end{equation}
The augmented Lagrangian function
is defined by
\[\label{augL}
{\cal L}^t_{\sigma}(x,\lambda):=F(x,\xi_t)+\displaystyle \frac{1}{2\sigma}\left[ \|\Pi_{\Re^p_+}(\lambda+\sigma G(x,\xi_t))\|^2-\|\lambda\|^2\right],\ \forall\, (x,\lambda)\in \Re^n\times \Re^p.
\]
Then  PMMSopt for Problem (\ref{eq:1}) may be described as follows.\mbox{}\\[4pt]
{\bf PMMSopt}: A  proximal method of multipliers  for solving Problem (\ref{eq:1}).
\begin{description}
\item[Step 0 ] Input $\lambda^0=0 \in \Re^p$ and $x^0 \in \Re^n$. Set $t:=0$.
\item[Step 1] Set
 \begin{equation}\label{xna}
\begin{array}{l}
x^{t+1}= \displaystyle\hbox{arg}\min \,\left\{ \L^t_{\sigma }(x,\lambda^t) +\displaystyle\frac{\alpha}{2}\|x-x^t\|^2,x \in X_0\right\}\\[5mm]
\lambda^{t+1}=[\lambda^t+\sigma G(x^{t+1},\xi_t)]_+.
\end{array}
\end{equation}
\item[Step 2] Set $t:=t+1$  and go to Step 1.
\end{description}
In the above algorithm, $[y]_+=\Pi_{\Re^p_+}[y]$ denotes  the projection of $y$ on to $\Re^p_+$ for any $y \in \Re^p$.
 Note that the iterations $x^t=x^t(\xi_{[t-1]})$ and $\lambda^t=\lambda^t(\xi_{[t-1]})$ are mappings of the history $\xi_{[t-1]}=(\xi_1,\ldots, \xi_{t-1})$ of the generated random process and hence are random.

As far as we are concerned, the main contributions of this paper can be summarized as follows.
\begin{itemize}

\item[(1)] When $\sigma=T^{-1/2}$, $\alpha=T^{1/2}$,  under mild assumptions, it is proved that the regret of objective function of $T$ iterations is  of order ${\rm O}(T^{-1/2})$, and the regret of constraint violation of $T$ iterations is  of order ${\rm O}(T^{-1/2})$.
    \item[(2)]
When $\sigma=T^{-1/2}$, $\alpha=T^{1/2}$,  under mild assumptions, it is proved that the following probability guarantees hold:
\begin{equation}\label{eq:CPG}
{\rm Pr} \left[\displaystyle
\frac{1}{T}\sum_{t=0}^{T-1}G_i(x^t,\xi_t) \leq \omega_c (T)\right]
\geq 1-e^{-T^{1/4}}
\end{equation}
 for all $i\in \{1,\ldots,p\}$, with
$$
\omega_c (T)=\left[8\left(1+2\sqrt{p}\kappa_g^2\right)\displaystyle
\frac{\nu_g^2}{\epsilon_0}\right] T^{-1/4}+{\rm o}(T^{-1/4});
$$
 and
\begin{equation}\label{eq:ObjPG}
{\rm Pr} \left[\displaystyle \displaystyle
\frac{1}{T}\sum_{t=0}^{T-1}F(x^t,\xi_t) \leq
\displaystyle\displaystyle \frac{1}{T} \sum_{t=0}^{T-1}F(\widehat
x,\xi_t) + \omega_o (T)\right] \geq 1-e^{-T^{1/4}},
\end{equation}
with
$$
\omega_o (T)=\displaystyle \frac{8\sqrt{2}\nu_g^3}{\epsilon_0} T^{-1/8}+{\rm o}(T^{-1/8}).
$$
\end{itemize}
The remaining parts of this paper are organized as follows. In Section 2, we develop properties of PMMSopt, which play a key role in the analysis for objective regret, constraint violation regret. In Section 3, we establish bounds of objective regret, constraint violation regret of PMMSopt for Problem (\ref{eq:1}). In Section 4, we develop probability guarantees for objective reduction and constraint violation of PMMSopt.  We draw a conclusion and give some discussions in Section \ref{Sec5}.


\section{Properties of PMMSopt}\label{sec-3}
\setcounter{equation}{0}
In this section, we develop properties of PMMSopt, which will be used in the analysis for objective regret and constraint violation regret.
\begin{lemma}\label{Lemma 2}
Let $(x^t,\lambda^t)$ be generated by PMMSopt and Assumption A(1)--Assumption (A3)  be satisfied. Then
\begin{equation}\label{eq:lambda2}
\|\lambda^{t+1}\|^2 \leq \|\lambda^t\|^2 +2 \sigma \langle \lambda^t, G(x^{t+1}, \xi_t)\rangle+\sigma^2 \nu_g^2,
\end{equation}
and
\begin{equation}\label{eq:2lambda}
\|\lambda^t\|- \sigma  \nu_g\leq \|\lambda^{t+1}\| \leq \|\lambda^t\|+ \sigma  \nu_g.
\end{equation}
\end{lemma}
{\bf Proof}. Noting that for any $a \in \Re$, $[a]_+^2 \leq a^2$, we have
$$
\begin{array}{ll}
\|\lambda^{t+1}\|^2 &= \displaystyle \sum_{i=1}^p [\lambda^t_i+\sigma G_i(x^{t+1},\xi_t)]_+^2\\[6pt]
&\leq \displaystyle \sum_{i=1}^p [\lambda^t_i+\sigma G_i(x^{t+1},\xi_t)]^2\\[6pt]
&= \displaystyle \sum_{i=1}^p\left( [\lambda^t_i]^2+2\lambda^t_i\sigma G_i(x^{t+1},\xi_t)+\sigma^2  G_i(x^{t+1},\xi_t)^2\right)\\[12pt]
&\leq \|\lambda^t\|^2 +2 \sigma \langle \lambda^t, G(x^{t+1}, \xi_t)\rangle+\sigma^2 \nu_g^2.
\end{array}
$$
It follows from the nonexpansion property of the projection $\Pi_{\Re^p_+}(\cdot)$, we have
$$
\|\lambda^{t+1}-\lambda^t\|=\|[\lambda^t+\sigma G(x^{t+1},\xi_t)]_+-[\lambda^t]_+\|\leq
\sigma \|G(x^{t+1},\xi_t)\|,
$$
which implies (\ref{eq:2lambda}).
The proof is completed. \hfill $\Box$
\begin{lemma}\label{Lemma 3}
Let $(x^t,\lambda^t)$ be generated by PMMSopt and Assumption (A1)--Assumption (A4)  be satisfied. Then, for  any stochastic subgradient $v_i(x^t,\xi_t)$   of $G_i(\cdot, \xi_t)$ at $x^t$, $i=1,\ldots,p$,
\begin{equation}\label{eq:cineq0a}
\begin{array}{ll}
\displaystyle \sum_{t=0}^{T-1} G_i(x^t,\xi_t)
&\leq
\displaystyle \frac{1}{\sigma} \lambda^{T}_i+ \sum_{t=0}^{T-1}\|v_i(x^t,\xi_t)\|\|x^{t+1}-x^t\|
\end{array}
\end{equation}
and
\begin{equation}\label{eq:cineqa}
\mathbb{E}\left[\displaystyle
\sum_{t=0}^{T-1}G_i(x^t,\xi_t)\right]\leq \displaystyle
\frac{1}{\sigma}\mathbb{E}[\lambda^{T}_i]+\kappa_g
\sum_{t=0}^{T-1}\mathbb{E}[\|x^{t+1}-x^t\|].
\end{equation}
\end{lemma}
{\bf Proof}. From the definition $\lambda^{t+1}_i=[\lambda^t_i+\sigma G_i(x^{t+1},\xi_t)]_+$, we have from  the convexity of $G_i(\cdot,\xi_t)$, for any stochastic subgradient $v_i(x^t,\xi_t)$ of $G_i(\cdot,\xi_t)$ at $x^t$ that
$$
\begin{array}{ll}
\lambda^{t+1}_i & \geq \lambda^t_i+\sigma G_i(x^{t+1},\xi_t)\\[10pt]
&\geq \lambda^t_i+\sigma( G_i(x^t,\xi_t)+\langle v_i(x^t,\xi_t), x^{t+1}-x^t\rangle )\\[10pt]
& \geq \lambda^t_i+\sigma(G_i(x^t,\xi_t)- \|v_i(x^t,\xi_t)\|\|x^{t+1}-x^t\|),
\end{array}
$$
which, from $\lambda^0=0$,  implies that
$$
\begin{array}{ll}
\displaystyle \sum_{t=0}^{T-1} G_i(x^t,\xi_t)&\leq  \displaystyle \frac{1}{\sigma} \lambda^{T}_i+ \sum_{t=0}^{T-1}\|v_i(x^t,\xi_t)\|\|x^{t+1}-x^t\|\\[10pt]
&\leq
\displaystyle \frac{1}{\sigma} \lambda^{T}_i+\kappa_g \sum_{t=0}^{T-1}\|x^{t+1}-x^t\|.
\end{array}
$$
Taking expectation operation in the both sides of the above inequality and using Assumption (A4), we obtain
(\ref{eq:cineqa}). \hfill $\Box$
\begin{lemma}\label{Lemma 4} Let $(x^t,\lambda^t)$ be generated by PMMSopt and Assumption (A1)--Assumption (A3) be satisfied. Then for any $x \in X_0$,
\begin{equation}\label{eq:l4}
\begin{array}{l}
F(x^{t+1},\xi_t)+\displaystyle \frac{1}{2\sigma}\|\lambda^{t+1}\|^2+\displaystyle \frac{\alpha}{2}\|x^{t+1}-x^t\|^2\\[8pt]
\leq F(x,\xi_t)+\displaystyle \frac{1}{2\sigma}\|[\lambda^t+\sigma
G(x,\xi_t)]_+\|^2+\displaystyle
\frac{\alpha}{2}[\|x-x^t\|^2-\|x-x^{t+1}\|^2].\end{array}
\end{equation}
\end{lemma}
{\bf Proof}. Since (\ref{eq:l4}) is equivalent to  the following inequality
\begin{equation}\label{eq:l4a}
{\cal L}^t_{\sigma}(x^{t+1},\xi_t)
+\displaystyle \frac{\alpha}{2}\|x^{t+1}-x^t\|^2
\leq  {\cal L}^t_{\sigma}(x,\xi_t)+\displaystyle \frac{\alpha}{2}[\|x-x^t\|^2-\|x-x^{t+1}\|^2],\,\, \forall x \in X_0,
\end{equation}
we only need to prove (\ref{eq:l4a}).
 Define
$$
\psi^t(x)=\delta_{X_0}(x)+{\cal L}^t_{\sigma}(x,\xi_t)+\displaystyle \frac{\alpha}{2}\|x-x^t\|^2
$$
 and
 $$
 q^t(x)={\cal L}^t_{\sigma}(x,\xi_t)-{\cal L}^t_{\sigma}(x^{t+1},\xi_t)-\langle \nabla_x {\cal L}^t_{\sigma}(x^{t+1},\xi_t),
 x-x^{t+1} \rangle.
 $$
Then
$$
\psi^t(x)=[{\cal L}^t_{\sigma}(x^{t+1},\xi_t)-\langle \nabla_x {\cal L}^t_{\sigma}(x^{t+1},\xi_t),
 x-x^{t+1}\rangle]+\displaystyle \frac{\alpha}{2}\|x-x^t\|^2+q^t(x)+\delta_{X_0} (x).
$$
Let
$$
q^t_0(x)=\delta_{X_0} (x)+q^t(x)+\langle \nabla_x {\cal L}^t_{\sigma}(x^{t+1})(x-x^{t+1})\rangle+\displaystyle \frac{1}{2}\|x-x^t\|^2-\displaystyle \frac{1}{2}\|x-x^{t+1}\|^2.
$$
Then $q^t_0$ is convex and
$$
\psi^t(x)={\cal L}^t_{\sigma}(x^{t+1},\xi_t)+\displaystyle \frac{1}{2}\|x-x^{t+1}\|^2+q^t_0(x).
$$
Since $0 \in \partial \psi^t(x^{t+1})$, we have $0 \in \partial q^t_0(x^{t+1})$ and $q^t_0$ arrives its minimum value at $x^{t+1}$.  Then for every $x \in \Re^n$,
$
q^t_0(x) \geq q^t_0(x^{t+1})= \delta_{X_0}(x^{t+1})+\displaystyle \frac{\alpha}{2}\|x^{t+1}-x^t\|^2
$, which is equivalent to
$$
\delta_{X_0}(x)+{\cal L}^t_{\sigma}(x,\xi_t)+\displaystyle \frac{\alpha}{2}\|x-x^t\|^2-\displaystyle \frac{\alpha}{2}\|x-x^{t+1}\|^2\geq \delta_{X_0}(x^{t+1})+{\cal L}^t_{\sigma}(x^{t+1},\xi_t)+\displaystyle \frac{\alpha}{2}\|x^{t+1}-x^t\|^2.
$$
The proof is completed. \hfill $\Box$\\

 In order to give a bound for $\displaystyle\sum_{t=0}^{T-1} G_i(x^t,\xi_t)$ in (\ref{eq:cineq0a}), we need to estimate an upper bound of $\displaystyle\sum_{t=0}^{T-1}\|x^{t+1}-x^t\|$, which is given in the following lemma.
 \begin{lemma}\label{Corollary 1}
 Let $(x^t,\lambda^t)$ be generated by PMMSopt and Assumption (A1)--Assumption (A4) be satisfied. Then for $2\alpha-p\kappa_g^2\sigma>0$,
\begin{equation}\label{eq:c1}
\|x^{t+1}-x^t\|\leq \displaystyle \frac{1}{2\alpha-p\kappa_g^2\sigma
}[2\kappa_f+\displaystyle \sqrt{p}\kappa_g\|\lambda^t\|+
\nu_g\sqrt{p}\kappa_g\sigma]
\end{equation}
\end{lemma}
{\bf Proof}. Let us denote
$$
B=\lambda^t+\sigma G(x^t,\xi_t), A=B+\sigma V(x^t,\xi_t)^T(x^{t+1}-x^t),
$$
where $V(x^t,\xi_t)=(v_1(x^t,\xi_t),\ldots,v_p(x^t,\xi_t))$.
Choosing $x=x^t$ in  Lemma \ref{Lemma 4}, we obtain
\begin{equation}\label{eq:dx}
\begin{array}{ll}
{\alpha} \|x^{t+1}-x^t\|^2 & \leq F(x^t,\xi_t)-F(x^{t+1},\xi_t)+\displaystyle \frac{1}{2\sigma}\|[\lambda^t+\sigma G(x^t,\xi_t)]_+\|^2\\[8pt]
& \quad -\displaystyle \frac{1}{2\sigma}\|[\lambda^t+\sigma
G(x^{t+1},\xi_t)]_+\|^2.
\end{array}
\end{equation}
Noting that
$$
 G(x^{t+1},\xi_t)\geq G(x^t,\xi_t)+V(x^t,\xi_t)^T(x^{t+1}-x^t),
$$
we have
$$
[\lambda^t+\sigma G(x^{t+1},\xi_t)]_+\geq
[\lambda^t+\sigma[G(x^t,\xi_t)+V(x^t,\xi_t)^T(x^{t+1}-x^t)]]_+=A_+
$$
and in turn,
\begin{equation}\label{dxa}
\|[\lambda^t+\sigma G(x^{t+1},\xi_t)]_+\|^2 \geq \|A_+\|^2.
\end{equation}
Therefore, we get from (\ref{eq:dx}) that
\begin{equation}\label{eq:dxb}
\begin{array}{ll}
{\alpha}\|x^{t+1}-x^t\|^2 & \leq F(x^t,\xi_t)-F(x^{t+1},\xi_t)+\displaystyle \frac{1}{2\sigma}\|[\lambda^t+\sigma G(x^t,\xi_t)]_+\|^2\\[6pt]
& \quad \quad -\displaystyle \frac{1}{2\sigma}\|[\lambda^t+\sigma G(x^{t+1},\xi_t)]_+\|^2\\[6pt]
&  \leq \langle v_0(x^t,\xi_t),x^t-x^{t+1}\rangle+\displaystyle \frac{1}{2\sigma}\|B_+\|^2-\displaystyle \frac{1}{2\sigma}\|A_+\|^2\\[6pt]
&  = \langle v_0(x^t,\xi_t),x^t-x^{t+1}\rangle+\displaystyle \frac{1}{2\sigma}[\|B_+\|+\|A_+\|][\|B_+\|-\|A_+\|]\\[6pt]
&  \leq \kappa_f\|x^t-x^{t+1}\|+\displaystyle \frac{1}{2\sigma}[\|B_+\|+\|A_+\|][\|B_+-A_+\|]\\[6pt]
&  \leq \kappa_f\|x^t-x^{t+1}\|+\displaystyle \frac{1}{2\sigma}[\|B\|+\|A\|][\|B-A\|]\\[6pt]
&  \leq \kappa_f\|x^t-x^{t+1}\|+\displaystyle \frac{1}{2\sigma}[\|\lambda^t\|+\sigma \|G(x^t,\xi_t)\|+\|B-A\|][\|B-A\|]\\[6pt]
&  \leq \kappa_f\|x^t-x^{t+1}\|+\displaystyle \frac{1}{2\sigma}[\|\lambda^t\|+\sigma \|G(x^t,\xi_t)\|\\[6pt]
& \quad \quad +\sigma \|V(x^t,\xi_t)^T(x^{t+1}-x^t)\|][\sigma\|V(x^t,\xi_t)^T(x^{t+1}-x^t)\|]\\[6pt]
&  \leq \kappa_f\|x^t-x^{t+1}\|+\displaystyle \frac{1}{2}\|\lambda^t\|\|V(x^t,\xi_t)^T(x^{t+1}-x^t)\|\\[6pt]
&\quad \quad +[\displaystyle \frac{1}{2}\|G(x^t,\xi^t)\|+ \displaystyle \frac{1}{2}\|V(x^t,\xi_t)^T(x^{t+1}-x^t)\|][\sigma\|V(x^t,\xi_t)^T(x^{t+1}-x^t)\|]\\[6pt]
&\leq \kappa_f\|x^t-x^{t+1}\|+\displaystyle \frac{1}{2}\sqrt{p}\kappa_g\|\lambda^t\|\|x^{t+1}-x^t\|\\[6pt]
&\quad \quad +[\displaystyle \frac{1}{2}\nu_g+ \displaystyle \frac{1}{2}\sqrt{p}\kappa_g\|x^{t+1}-x^t\|][\sigma\sqrt{p}\kappa_g\|x^{t+1}-x^t\|],\\[6pt]
\end{array}
\end{equation}
which implies
$$
\|x^{t+1}-x^t\|\leq \displaystyle \frac{1}{2\alpha-p\kappa_g^2\sigma
}[2\kappa_f+\displaystyle \sqrt{p}\kappa_g\|\lambda^t\|+
\nu_g\sqrt{p}\kappa_g\sigma]
$$
The proof is completed. \hfill $\Box$
\begin{corollary}\label{Corollary 2}
 Let $(x^t,\lambda^t)$ be generated by PMMSopt and Assumption (A1)--Assumption (A4) be satisfied. Then for $2\alpha-p\kappa_g^2\sigma>0$,
\begin{equation}\label{eq:c2}
\begin{array}{ll}
\displaystyle \sum_{t=0}^{T-1} G_i(x^t,\xi_t)
&\leq
\displaystyle \frac{1}{\sigma} \lambda^{T}_i+  \displaystyle \frac{2\kappa_g\kappa_f}{2\alpha-p\kappa_g^2\sigma }T\\[10pt]
& \quad + \displaystyle
\frac{\sqrt{p}\kappa_g^2}{2\alpha-p\kappa_g^2\sigma }\displaystyle
\sum_{t=0}^{T-1}\|\lambda^t\| +\displaystyle
\frac{\sqrt{p}\nu_g\kappa_g^2\sigma}{2\alpha-p\kappa_g^2\sigma }T.
\end{array}
\end{equation}
\end{corollary}
\begin{lemma}\label{Lemma 6}
Let $(x^t,\lambda^t)$ be generated by PMMSopt and Assumption (A1), Assumption (A2) and Assumption  (A5) be satisfied. Then for any $t_2 \leq t_1-1$ where $t_1$ and $t_2$ are positive integers,
\begin{equation}\label{eq:slaterCondlamd}
\displaystyle \mathbb{E}\left[\langle \lambda^{t_1},G(\widehat x, \xi_{t_1}) \rangle \,|\, \xi_{[t_2]}\right]
\leq -\epsilon_0 \mathbb{E}\left[\|\lambda^{t_1}\| \,|\, \xi_{[t_2]}\right].
\end{equation}
\end{lemma}
{\bf Proof}. To prove this lemma, we first show that
$$
\mathbb{E}\left[\lambda^{t_1}_iG_i(\widehat x, \xi_{t_1}) \,|\, \xi_{[t_2]}\right]
\leq -\epsilon_0 \mathbb{E}\left[\lambda^{t_1}_i \,|\, \xi_{[t_2]}\right].
$$
For $i \in \{1,\ldots,p\}$, note that
$\lambda^{t_1}_i=\lambda^{t_1}_i(\xi_{[t_1-1]})$ and $G_i(\widehat
x,\xi_{t_1})$ is independent of $\xi_{[t_1-1]}$ and
$\xi_{[t_2]}\subseteq \xi_{[t_1-1]}$ for $t_2 \leq t_1-1$. We have
$$
\begin{array}{ll}
\mathbb{E}\left[\lambda^{t_1}_iG_i(\widehat x, \xi_{t_1}) \,|\, \xi_{[t_2]}\right] &=\mathbb{E} \left\{\mathbb{E}\left[\lambda^{t_1}_iG_i(\widehat x, \xi_{t_1}) \,|\, \xi_{[t_1-1]}\right]\,|\,\xi_{[t_2]}\right\}\\[10pt]
&=\mathbb{E}\left\{\lambda^{t_1}_i\mathbb{E}\left[G_i(\widehat x, \xi_{t_1}) \right]\,|\,\xi_{[t_2]}\right\}\\[10pt]
&=\mathbb{E}\left[G_i(\widehat x, \xi_{t_1})\right]\left[\lambda^{t_1}_i \,|\,\xi_{[t_2]}\right]\\[10pt]
&\leq -\epsilon_0 \mathbb{E}\left[\lambda^{t_1}_i \,|\, \xi_{[t_2]}\right].
\end{array}
$$
Making a sum over $i \in \{1,\ldots,p\}$ yields
$$
\displaystyle \mathbb{E}\left[\langle \lambda^{t_1},G(\widehat x, \xi_{t_1}) \rangle \,|\, \xi_{[t_2]}\right]
\leq
-\epsilon_0 \mathbb{E} \left[ \sum_{i=1}^p \lambda^{t_1}_i \,|\, \xi_{[t_2]}\right]
\leq -\epsilon_0 \mathbb{E}\left[\|\lambda^{t_1}\| \,|\, \xi_{[t_2]}\right].
$$
The proof is completed.\hfill $\Box$\\
\begin{lemma}\label{Lemma 7}
 Let $s > 0$ be an arbitrary integer. Let Assumption (A1)--Assumption (A4) be satisfied. At each round  $t \in \{1,2,\ldots\}$ in  PMMSopt. For
 \begin{equation}\label{eq:theta9}
 \vartheta (\sigma,\alpha,s)=\displaystyle \frac{\epsilon_0\sigma s}{2}+ \nu_g\sigma(s-1)+\displaystyle \frac{\alpha D_0^2}{\epsilon_0s}+\displaystyle \frac{2\nu_f}{\epsilon_0}+
 \displaystyle \frac{\sigma \nu_g^2}{\epsilon_0},
 \end{equation}
 the following  holds
\begin{equation}\label{eq:6}
|\|\lambda^{t+1}\|-\|\lambda^t\||\leq \sigma \nu_g
\end{equation}
and
\begin{equation}\label{eq:7}
\mathbb{E}\left [ \|\lambda^{t+s}\|-\|\lambda^t\| \,|\, \xi_{[t-1]}\right]
\leq \left
\{
\begin{array}{ll}
s \sigma \nu_g & \mbox{if } \|\lambda^t\| < \vartheta (\sigma,\alpha,s),\\[6pt]
-s \displaystyle \frac{\sigma \epsilon_0}{2} & \mbox{if } \|\lambda^t\| \geq  \vartheta (\sigma,\alpha,s).
\end{array}
\right.
\end{equation}
\end{lemma}
{\bf Proof}. Inequality (\ref{eq:6}) follows from Lemma \ref{Lemma 2}. We only need to establish (\ref{eq:7}).
Since it is obvious that
$$
 \mathbb{E}\left[\|\lambda^{t+s}\|-\|\lambda^t\|\,|\,\xi_{[t-1]}\right]\leq s \sigma \nu_g
$$
 when $\|\lambda^t\|< \vartheta (\sigma,\alpha,s)$, it remains to prove
$$ \mathbb{E}\left [ \|\lambda^{t+s}\|-\|\lambda^t\| \,|\, \xi_{[t-1]}\right]
\leq -s \displaystyle \frac{\sigma \epsilon_0}{2}
$$
when $\|\lambda^t\| \geq  \vartheta (\sigma,\alpha,s)$.

For given positive integer $s$, suppose $\|\lambda^t\| \geq  \vartheta (\sigma,\alpha,s)$. For any $l \in \{t,t+1,\ldots, t+s-1\}$, one has
$$
\begin{array}{l}
F(x^{l+1},\xi_l)+\displaystyle \frac{1}{2\sigma}\|\lambda^{l+1}\|^2+\displaystyle \frac{\alpha}{2}\|x^{l+1}-x^l\|^2\\[10pt]
\quad \,\leq F(\widehat x, \xi_l)+\displaystyle \frac{1}{2\sigma}\|[\lambda^{l}+\sigma G(\widehat x,\xi_l)]_{+}\|^2+\displaystyle \frac{\alpha}{2}\left[\|\widehat x-x^l\|^2-\|\widehat x-x^{l+1}\|^2\right].
\end{array}
$$
Using Assumption (A3) and the following inequality
$$
\|[\lambda^{l}+\sigma G(\widehat x,\xi_l)]_{+}\|^2 \leq \|\lambda^l\|^2+2\sigma \langle \lambda^l,G(\widehat x,\xi_l)\rangle
+\sigma^2\|G(\widehat x,\xi_l)\|^2,
$$
we obtain
\begin{equation}\label{eq:8}
\begin{array}{ll}
\displaystyle \frac{1}{2\sigma} \left[\|\lambda^{l+1}\|^2-\|\lambda^l\|^2\right]
&\leq \left(F(\widehat x,\xi_l)-F(x^{l+1},\xi_l)\right)\\[10pt]
& \quad \quad +\displaystyle \frac{1}{2\sigma}\left[\|[\lambda^{l}+\sigma G(\widehat x,\xi_l)]_{+}\|^2-\|\lambda^l\|^2\right]\\[10pt]
&\quad \quad -\displaystyle \frac{\alpha}{2}\|x^{l+1}-x^l\|^2+\displaystyle \frac{\alpha}{2}\left[\|\widehat x-x^l\|^2-\|\widehat x-x^{l+1}\|^2\right]\\[10pt]
&\leq \nu_f+ \langle \lambda^l,G(\widehat x,\xi_l)\rangle
+\displaystyle \frac{\sigma}{2}\|G(\widehat x,\xi_l)\|^2\\[10pt]
&\quad \quad -\displaystyle \frac{\alpha}{2}\|x^{l+1}-x^l\|^2+\displaystyle \frac{\alpha}{2}\left[\|\widehat x-x^l\|^2-\|\widehat x-x^{l+1}\|^2\right].
\end{array}
\end{equation}
Making  a summation of (\ref{eq:8}) over $\{t,t+1, t+s-1\}$ and taking conditional expectation on $\xi_{[t-1]}$ , we obtain from
Lemma \ref{Lemma 7} that
\begin{equation}\label{eq:9}
\begin{array}{l}
\displaystyle \frac{1}{2\sigma} \mathbb{E}\left[\|\lambda^{t+s}\|^2-\|\lambda^t\|^2\,|\, \xi_{[t-1]}\right]\\[10pt]
\leq \nu_f s + \displaystyle \frac{\sigma}{2}\nu_g^2s +\displaystyle \sum_{l=t}^{t+s-1} \mathbb{E}\left[\langle \lambda^l,G(\widehat x,\xi_l)\rangle\,|\, \xi_{[t-1]}\right]
\\[10pt]
\quad \quad +\displaystyle \frac{\alpha}{2}\mathbb{E}\left[\left(\|\widehat x-x^t\|^2-\|\widehat x-x^{t+s}\|^2\right)\,|\, \xi_{[t-1]}\right]\\[10pt]
\leq \nu_f s + \displaystyle \frac{\sigma}{2}\nu_g^2s -\epsilon_0\displaystyle \sum_{l=0}^{s-1} \mathbb{E}\left[\|\lambda^{t+l}\|\,|\, \xi_{[t-1]}\right]
\\[10pt]
\quad \quad +\displaystyle \frac{\alpha}{2}\mathbb{E}\left[\left(\|\widehat x-x^t\|^2-\|\widehat x-x^{t+s}\|^2\right)\,|\, \xi_{[t-1]}\right]\\[10pt]
\leq \nu_f s + \displaystyle \frac{\sigma}{2}\nu_g^2s -\epsilon_0\displaystyle \sum_{l=0}^{s-1} \mathbb{E}\left[\|\lambda^{t}\|-\sigma\nu_gl \,|\, \xi_{[t-1]}\right]
\\[10pt]
\quad \quad +\displaystyle \frac{\alpha}{2}\mathbb{E}\left[\left(\|\widehat x-x^t\|^2-\|\widehat x-x^{t+s}\|^2\right)\,|\, \xi_{[t-1]}\right]\\[10pt]
 \quad \quad (\mbox{from } \|\lambda^{t+1}\|\geq \|\lambda^t\|-\sigma \nu_g)\\[8pt]
\leq \nu_f s + \displaystyle \frac{\sigma}{2}\nu_g^2s +\displaystyle \frac{\alpha}{2}\mathbb{E}\left[\left(\|\widehat x-x^t\|^2-\|\widehat x-x^{t+s}\|^2\right)\,|\, \xi_{[t-1]}\right]\\[10pt]
 \quad \quad + \epsilon_0\sigma \nu_g \displaystyle \frac{s(s-1)}{2}
-\epsilon_0\displaystyle \sum_{l=0}^{s-1} \mathbb{E}\left[\|\lambda^{t}\| \,|\, \xi_{[t-1]}\right]
\end{array}
\end{equation}
From (\ref{eq:9}), we get from Assumption (A3) that
\begin{equation}\label{eq:10}
\begin{array}{l}
\mathbb{E}\left[\|\lambda^{t+s}\|^2\,|\, \xi_{[t-1]}\right]\leq
\mathbb{E}\left[\|\lambda^t\|^2\,|\, \xi_{[t-1]}\right]\\[10pt]
\quad \quad +2 \sigma \nu_f s+\sigma^2\nu_g^2s+\alpha \sigma D_0^2
+\epsilon_0\sigma^2\nu_gs(s-1)-2\epsilon_0\sigma s \mathbb{E}\left[\|\lambda^{t}\| \,|\, \xi_{[t-1]}\right]\\[10pt]
=\mathbb{E}\left[(\|\lambda^t\|-\displaystyle \frac{\epsilon_0\sigma}{2}s)^2\,|\, \xi_{[t-1]}\right]
-\displaystyle \frac{\epsilon_0^2\sigma^2}{4}s^2+\epsilon_0\sigma^2 \nu_gs(s-1)\\[10pt]
\quad \quad + \alpha \sigma D_0^2+2\sigma \nu_fs+\sigma^2\nu_g^2s-\epsilon_0\sigma s \mathbb{E}\left[\|\lambda^{t}\| \,|\, \xi_{[t-1]}\right]\\[10pt]
\leq \mathbb{E}\left[(\|\lambda^t\|-\displaystyle \frac{\epsilon_0\sigma}{2}s)^2\,|\, \xi_{[t-1]}\right]
-\displaystyle \frac{3\epsilon_0^2\sigma^2}{4}s^2\\[10pt]
\quad \quad +\left[\epsilon_0\sigma^2 \nu_gs(s-1)+\displaystyle \frac{\epsilon_0^2\sigma^2}{2}s^2+\alpha \sigma D_0^2+2\sigma \nu_fs+\sigma^2\nu_g^2s-\epsilon_0\sigma s \vartheta (\sigma,\alpha,s)\right]\\[10pt]
= \mathbb{E}\left[(\|\lambda^t\|-\displaystyle \frac{\epsilon_0\sigma}{2}s)^2\,|\, \xi_{[t-1]}\right]-\displaystyle \frac{3\epsilon_0^2\sigma^2}{4}s^2\\[10pt]
\leq \mathbb{E}\left[(\|\lambda^t\|-\displaystyle \frac{\epsilon_0\sigma}{2}s)^2\,|\, \xi_{[t-1]}\right].
\end{array}
\end{equation}
This implies that
$$
\mathbb{E}\left[\|\lambda^{t+s}\|\,|\, \xi_{[t-1]}\right]\leq
\mathbb{E}\left[\|\lambda^t\|\,|\, \xi_{[t-1]}\right]-\displaystyle \frac{\epsilon_0\sigma}{2}s.
$$
The proof is completed. \hfill $\Box$\\

The following lemmas come from Yu et.al \cite{YMNeely2017}, which
can be used to deal with the random process $\{\|\lambda^t\|\}$ and
probability analysis for objective regret and constraint violation
regret, respectively.
\begin{lemma}\label{Lemma 5}
Let $\{Z(t), t \geq 0\}$ be a discrete time stochastic process adapted to a filtration $\{{\cal F}(t), t\geq
0\}$ with $Z(0) = 0$ and ${\cal F}(0) = \{\emptyset, \Omega\}$. Suppose there exists an integer $t_0 >0$, real constants $\theta>0$, $\delta_{\max}>0$ and $ 0 <\zeta \leq \delta_{\max}$ such that
\begin{equation}\label{eq:Y9}
\begin{array}{rl}
|Z(t+1)-Z(t)| & \leq \delta_{\max}\mbox{ and }\\[12pt]
\mathbb{E}[Z(t+t_0)-Z(t)\,|\, {\cal F}(t)] & \leq \left
\{
\begin{array}{ll}
t_0 \delta_{\max} & \mbox{if } Z(t) < \theta\\[6pt]
-t_0\zeta & \mbox{if } Z(t) \geq \theta
\end{array}
\right.
\end{array}
\end{equation}
hold for all $t \in \{1,2,\ldots\}.$ Then the following properties are satisfied.
\begin{itemize}
\item[1.] The  following inequality holds
$$\mathbb{E}[Z(t)] \leq \theta +t_0 \delta_{\max}+t_0 \displaystyle \frac{4 \delta_{\max}^2}{\zeta}\log \left[ \displaystyle \frac{8 \delta_{\max}^2}{\zeta^2} \right], \forall t \in \{1,2,\ldots\}.$$
    \item[2.] For any constant $0 < \mu <1$, we have
    $$
    {\rm Pr}\left\{Z(t)\geq z\right\} \leq \mu, \forall t \in \{1,2,\ldots\},
    $$
    where
    $$
    z=\theta +t_0 \delta_{\max}+t_0 \displaystyle \frac{4 \delta_{\max}^2}{\zeta}\log \left[ \displaystyle \frac{8 \delta_{\max}^2}{\zeta^2}\right]+t_0 \displaystyle \frac{4 \delta_{\max}^2}{\zeta}\log\left(\displaystyle \frac{1}{\mu} \right).
    $$
\end{itemize}
\end{lemma}
\begin{lemma}\label{Lemma 9}
Let $\{Z(t),t\geq 0\}$ be a supermartingale adapted to a filtration $\{{\cal F}(t),t\geq 0\}$ with $Z(0)=0$ and ${\cal F}(0)=\{\emptyset, \Omega\}$, i.e. $\mathbb E[Z(t+1)\,|\, {\cal F}(t)]\leq Z(t),\forall t \geq 0$. Suppose there exists a constant $c>0$ such that $\{|Z(t+1)-Z(t)|>c\}\subseteq \{Y(t)>0\},  \forall t\geq 0$, where $Y(t)$ is process with $Y(t)$ adapted to ${\cal F}(t)$ for all $t\geq 0$.  Then, for all $z>0$, we have
$$
{\rm Pr}[Z(t)\geq z] \leq e^{-z^2/(2tc^2)}+\displaystyle \sum_{j=0}^{t-1} {\rm Pr}[Y(j)>0], \forall t \geq 1.
$$
\end{lemma}
\section{Regret analysis of PMMSopt}\label{sec3}
\setcounter{equation}{0}

In order to use Lemma \ref{Lemma 5} and Lemma \ref{Lemma 7} to
analyze regrets of PMMSopt for Problem (\ref{eq:1}), we introduce
the following notations. For $\theta=\vartheta (\sigma,\alpha,s)$,
$\delta_{\max}=\sigma \nu_g$ and $\zeta =\displaystyle
\frac{\sigma}{2}\epsilon_0$, and $t_0=s$, define
$$
\psi(\sigma,\alpha,s)=\theta +t_0 \delta_{\max}+t_0 \displaystyle \frac{4 \delta_{\max}^2}{\zeta}\log \left[ \displaystyle \frac{8 \delta_{\max}^2}{\zeta^2} \right]
$$
 and
$$
\phi (\sigma,\alpha,s,\mu)=\psi(\sigma,\alpha,s)+8  \displaystyle \frac{ \nu_g^2}{\epsilon_0} \log \left( \displaystyle \frac{1}{\mu}\right)\sigma s.
$$
 Then $\psi(\sigma,\alpha,s)$ is expressed as
$$
\begin{array}{ll}
\psi(\sigma,\alpha,s)= &
\vartheta (\sigma,\alpha,s)+  \left[\nu_g+\displaystyle \frac{8\nu_g^2}{\epsilon_0}\log\displaystyle \frac{32\nu_g^2}{\epsilon_0^2}\right] \sigma s\\[12pt]
&=\kappa_0+\kappa_1\displaystyle \frac{\alpha}{s}+ \kappa_2 s+\kappa_3
 \sigma+\kappa_4 \sigma s
 \end{array}
$$
 and
 $\phi (\sigma,\alpha,s,\mu)$ is expressed as
 $$
 \phi (\sigma,\alpha,s,\mu)=\kappa_0+\kappa_1\displaystyle \frac{\alpha}{s}+ \kappa_2 s+\kappa_3
 \sigma+\kappa_4 \sigma s+8 \displaystyle \frac{ \nu_g^2}{\epsilon_0} \log \left( \displaystyle \frac{1}{\mu}\right)\sigma s
 $$
 where
\begin{equation}\label{eq:notations}
\begin{array}{l}
\kappa_0=\displaystyle \frac{2\nu_f}{\epsilon_0},\,\,
\kappa_1=\displaystyle \frac{D_0^2}{\epsilon_0},\,\,
\kappa_2=0,\\[10pt]
\kappa_3=\displaystyle \frac{ \nu_g^2}{\epsilon_0}-\nu_g,\,\,
\kappa_4=\left[2\nu_g +\displaystyle \frac{\epsilon_0}{2}+\displaystyle \frac{8\nu_g^2}{\epsilon_0}\log \displaystyle \frac{32\nu_g^2}{\epsilon_0^2}\right].
\end{array}
\end{equation}

  Lemma \ref{Lemma 7} allows us to apply Lemma \ref{Lemma 5} to random process $Z(t) =\|\lambda^t\|$
and obtain $\mathbb{E}[\|\lambda^t\|] =\mbox{O}(1)$, $\forall t$ by taking $s=\lceil\sqrt{T}\rceil$,$\alpha= \sqrt{T}$, where $\lceil\sqrt{T}\rceil$
represents the smallest integer no less than $\sqrt{T}$. By Corollary \ref{Corollary 2}, this further implies the expected
constraint violation bound
$$
\mathbb{E}\left[\displaystyle \sum_{t=0}^{T-1}G_i(x^t,\xi_t)\right]\leq\mbox{O}(\sqrt{T}),i=1,\ldots,p
$$
 as summarized in the next theorem.
 \begin{theorem}\label{Theorem 1}
 (Expected Constraint Violation Bound). If $\alpha = \sqrt{T}$  and  $\sigma=T^{-1/2}$ in PMMSopt, then for
all $T \geq 1$, we have
\begin{equation}\label{eq:Cregret}
\mathbb{E}\left[\displaystyle \sum_{t=0}^{T-1}G_i(x^t,\xi_t)\right] \leq \kappa_c\sqrt{T},i=1,\ldots,p,
\end{equation}
where
$$
\kappa_c= \kappa_*+4\kappa_g\kappa_f +
2\sqrt{p}(\kappa_*+\nu_g)\kappa_g^2,
$$
with $\kappa_*=\kappa_0+\kappa_1+ \kappa_3
 +\kappa_4$.
\end{theorem}
{\bf Proof}. Define random process $Z(t)$ with $Z(0)=0$ and
$Z(t)=\|\lambda^t\|$,$t\geq 1$ and filtration ${\cal F}(t)$ with
${\cal F}(0)=\{\emptyset, \Omega\}$ and ${\cal F}(t)=\xi_{[t-1]}$,$t
\geq 1$. Then we have $Z(t)$ is adapted to ${\cal F}(t)$. It follows
from Lemma \ref{Lemma 7} that $Z(t)=\|\lambda^t\|$ satisfies Lemma
\ref{Lemma 5} with $\delta_{\max}=\sigma \nu_g$, $t_0=s$ and
$\zeta=\sigma \epsilon_0/2$. We have from Lemma \ref{Lemma 5} that
the following inequality holds for every $t$:
$$
\mathbb{E}\|\lambda^t\|\leq \psi(\sigma,\alpha,s)=\kappa_0+\kappa_1\displaystyle \frac{\alpha}{s}+ \kappa_3
 \sigma+\kappa_4 \sigma s.
$$
 Taking $s=\lceil\sqrt{T}\rceil$,$\alpha= \sqrt{T}$, $2\sqrt{T}-p\kappa_g^2\sigma> \sqrt{T}/2$ and $\sigma=T^{-1/2}$, we have $\mathbb{E}\|\lambda^t\|\leq \kappa_*$ for $t\in \{1,2,\ldots, T\}$, where
 $
 \kappa_*=\kappa_0+\kappa_1+ \kappa_3
 +\kappa_4$.

 From Corollary \ref{Corollary 2}, we have
 $$
\begin{array}{ll}
\mathbb{E}\left[\displaystyle \sum_{t=0}^{T-1}G_i(x^t,\xi_t)\right]
&\leq \displaystyle \frac{1}{\sigma} \mathbb{E}\lambda^{T}_i+
\displaystyle \frac{2\kappa_g\kappa_f}{2\alpha-p\kappa_g^2\sigma }T
+ \displaystyle \frac{\sqrt{p}\kappa_g^2}{2\alpha-p\kappa_g^2\sigma
}\displaystyle \sum_{t=0}^{T-1}\mathbb{E}\|\lambda^t\|
+\displaystyle \frac{\sqrt{p}\nu_g\kappa_g^2\sigma}{2\alpha-p\kappa_g^2\sigma }T\\[16pt]
& \leq \kappa_*\sqrt{T}+4\kappa_g\kappa_f\sqrt{T} +
2\sqrt{p}\kappa_*\kappa_g^2\sqrt{T}+2\sqrt{p}\nu_g\kappa_g^2\leq
\kappa_c\sqrt{T},
\end{array}
$$
which implies (\ref{eq:Cregret}). The proof is completed. \hfill $\Box$\\
\begin{lemma}\label{Lemma 8}
 Let $(x^t,\lambda^t)$ be generated by PMMSopt and Assumption (A1)--Assumption (A4) be satisfied. Then for $x \in X_0$,
\begin{equation}\label{eq:n13}
\begin{array}{ll}
\displaystyle \sum_{t=0}^{T-1} F(x^t,\xi_t) & \leq \displaystyle \sum_{t=0}^{T-1} F(x,\xi_t) +\displaystyle \frac{\kappa_f^2}{2\alpha}T +\displaystyle \frac{\alpha}{2}D_0^2+\displaystyle \frac{\sigma}{2}\nu_g^2T\\[12pt]
&\quad  + \displaystyle\displaystyle \sum_{t=0}^{T-1}\left[ \langle \lambda^t,G(x,\xi_t)\rangle\right].
\end{array}
\end{equation}
\end{lemma}
{\bf Proof}. In view of  (\ref{eq:l4}), we have for $x\in X_0$,
\begin{equation}\label{eq:h1}
\begin{array}{ll}
F(x^t,\xi_t) & \leq F(x,\xi_t)+[F(x^t,\xi_t)-F(x^{t+1},\xi_t)]-\displaystyle \frac{1}{2\sigma}\left[\|\lambda^{t+1}\|^2-\|\lambda^t\|^2\right]\\[6pt]
& \quad \, + \displaystyle \frac{1}{2\sigma}\left[\|[\lambda^{t}+\sigma G(x,\xi_t)]_+\|^2-\|\lambda^t\|^2\right]\\[6pt]
& \quad +\displaystyle \frac{\alpha}{2}\left[\|x-x^t\|^2-\|x-x^{t+1}\|^2-\|x^{t+1}-x^t\|^2\right]\\[6pt]
& = F(x,\xi_t)+\left[F(x^t,\xi_t)-F(x^{t+1},\xi_t)-\displaystyle \frac{\alpha}{2}\|x^{t+1}-x^t\|^2\right
]\\[6pt]
&\quad -\displaystyle \frac{1}{2\sigma}\left[\|\lambda^{t+1}\|^2-\|\lambda^t\|^2\right] + \displaystyle \frac{1}{2\sigma}\left[\|[\lambda^{t}+\sigma G(x,\xi_t)]_+\|^2-\|\lambda^t\|^2\right]\\[12pt]
& \quad +\displaystyle \frac{\alpha}{2}\left[\|x-x^t\|^2-\|x-x^{t+1}\|^2\right].
\end{array}
\end{equation}
From the convexity of $F(\cdot,\xi_t)$, we obtain from Assumption (A4) that
\begin{equation}\label{eq:h2}
\begin{array}{l}
\left[F(x^t,\xi_t)-F(x^{t+1},\xi_t)-\displaystyle \frac{\alpha}{2}\|x^{t+1}-x^t\|^2\right]\\[6pt]
\leq \left[\langle v_0(x^t,\xi_t),x^t-x^{t+1}\rangle-\displaystyle \frac{\alpha}{2}\|x^{t+1}-x^t\|^2\right]\\[6pt]
=-\displaystyle \frac{\alpha}{2}\left\{\|x^{t+1}-x^t+v_0(x^t,\xi_t)/\alpha\|^2-\|v_0(x^t,\xi_t)\|^2/\alpha^2\right\}\\[6pt]
\leq \displaystyle \frac{\kappa_f^2}{2\alpha}.
\end{array}
\end{equation}
Since $[a]_+^2\leq a^2$ for scalar $a \in \Re$, we obtain
\begin{equation}\label{eq:h3}
\begin{array}{l}
\|[\lambda^t+\sigma G(x,\xi_t)]_+\|^2-\|\lambda^t\|^2\\[6pt]
\leq \|\lambda^t+\sigma G(x,\xi_t)\|^2-\|\lambda^t\|^2\\[6pt]
=2\sigma \langle \lambda^t,G(x,\xi_t)\rangle+\sigma^2\|G(x,\xi_t)\|^2.
\end{array}
\end{equation}
Substituting (\ref{eq:h2}) and (\ref{eq:h3}) into (\ref{eq:h1}), we get
\begin{equation}\label{eq:h4}
\begin{array}{ll}
F(x^t,\xi_t) &
\leq F(x,\xi_t)+\displaystyle \frac{\kappa_f^2}{2\alpha}-\displaystyle \frac{1}{2\sigma}\left[\|\lambda^{t+1}\|^2-\|\lambda^t\|^2\right]\\[6pt]
&\quad  + \displaystyle \frac{1}{2}\left[2 \langle \lambda^t,G(x,\xi_t)\rangle+\sigma\|G(x,\xi_t)\|^2\right]\\[6pt]
& \quad +\displaystyle \frac{\alpha}{2}\left[\|x-x^t\|^2-\|x-x^{t+1}\|^2\right]\\[6pt]
\end{array}
\end{equation}
Making a summation, we obtain from Assumption (A3)
\begin{equation}\label{eq:h4}
\begin{array}{ll}
\displaystyle \sum_{t=0}^{T-1}F(x^t,\xi_t) &
\leq \displaystyle \sum_{t=0}^{T-1} F(x,\xi_t)+\displaystyle \frac{\kappa_f^2}{2\alpha}T-\displaystyle \frac{1}{2\sigma}\left[\|\lambda^{T}\|^2-\|\lambda^0\|^2\right]\\[6pt]
&\quad  + \displaystyle \frac{1}{2}\displaystyle \sum_{t=0}^{T-1}\left[2 \langle \lambda^t,G(x,\xi_t)\rangle+\sigma\|G(x,\xi_t)\|^2\right]\\[6pt]
& \quad +\displaystyle \frac{\alpha}{2}\left[\|x-x^0\|^2-\|x-x^{T}\|^2\right]\\[6pt]
&\leq \displaystyle \sum_{t=0}^{T-1} F(x,\xi_t)+\displaystyle \frac{\kappa_f^2}{2\alpha}T +\displaystyle \frac{\alpha}{2}D_0^2+\displaystyle \frac{\sigma}{2}\nu_g^2T\\[12pt]
&\quad  + \displaystyle\displaystyle \sum_{t=0}^{T-1}\left[ \langle \lambda^t,G(x,\xi_t)\rangle\right],
\end{array}
\end{equation}
which follows from $\lambda^0=0$. \hfill $\Box$\\
\begin{theorem}\label{Theorem 2}
 (Expected Regret Bound) Let $x^*\in \Phi$ be any fixed solution that satisfies
$$
 x^* = \mbox{argmin}_{x \in \Phi}\displaystyle \sum_{t=0}^{T-1} f_t(x).
 $$
 If $\alpha=\sqrt{T}$,$\sigma=T^{-1/2}$ in PMMSopt and Assumption (A1)--Assumption (A4) are satisfied. Then for all $T \geq 2$,
 \begin{equation}\label{eq:t2}
 \mathbb{E}\left[\displaystyle \sum_{t=0}^{T-1} F(x^t,\xi_t) \right]\leq \mathbb{E}\left[\displaystyle \sum_{t=0}^{T-1} F(x^*,\xi_t) \right]+\kappa_o\sqrt{T},
  \end{equation}
  where
  $$\kappa_o=\displaystyle \frac{\kappa_f^2}{2} +\displaystyle \frac{1}{2}D_0^2+\displaystyle \frac{1}{2}\nu_g^2.
  $$
  \end{theorem}
  {\bf Proof}. For fixed $T \geq 2$, taking $x =x^*$ in (\ref{eq:n13}), we get
\begin{equation}\label{eq:n13a}
\begin{array}{ll}
\displaystyle \sum_{t=0}^{T-1} \mathbb{E} F(x^t,\xi_t) & \leq \displaystyle \sum_{t=0}^{T-1}\mathbb{E} F(x^*,\xi_t) +\displaystyle \frac{\kappa_f^2}{2\alpha}T +\displaystyle \frac{\alpha}{2}D_0^2+\displaystyle \frac{\sigma}{2}\nu_g^2T\\[12pt]
&\quad  + \displaystyle\displaystyle \sum_{t=0}^{T-1}\mathbb{E}\left[ \langle \lambda^t,G(x^*,\xi_t)\rangle\right].
\end{array}
\end{equation}
Since $G(x^*,\xi_t)$ is independent on $\lambda^t$, which is determined by $\xi_{[t-1]}$, we have
$$
\mathbb{E}\left[ \langle \lambda^t,G(x^*,\xi_t)\rangle\right]=\mathbb{E}\left[ \langle \lambda^t,\mathbb{E}[G(x^*,\xi_t)|\xi_{[t-1]}]\rangle\right]\leq 0,
$$
which follows from the fact that $\lambda^t \geq 0$ and $\mathbb{E}[G(x^*,\xi_t)|\xi_{[t-1]}]=\mathbb{E}G(x^*,\xi_t)\leq0$.
 Thus we have from (\ref{eq:n13a}) and $\alpha=\sqrt{T}$,$\sigma=T^{-1/2}$ that
 $$
 \displaystyle \sum_{t=0}^{T-1} \mathbb{E} F(x^t,\xi_t) \leq \displaystyle \sum_{t=0}^{T-1}\mathbb{E} F(x^*,\xi_t) +\left[\displaystyle \frac{\kappa_f^2}{2} +\displaystyle \frac{1}{2}D_0^2+\displaystyle \frac{1}{2}\nu_g^2\right]\sqrt{T}.
 $$
   The proof is completed. \hfill $\Box$

\section{High probability performance analysis}
\setcounter{equation}{0}
\quad \, First of all, we will use (\ref{eq:c2}) and part 2 of  Lemma \ref{Lemma 7} to establish a high probability constraint violation bound.
\begin{theorem}\label{th:constr3} Let $\eta\in (0,1)$.  Let $(x^t,\lambda^t)$ be generated by PMMSopt, and Assumption (A1)--Assumption (A4) be satisfied. If $\sigma=T^{-1/2},\, \alpha=T^{1/2}$ in PMMSopt, then for all $i\in \{1,\ldots,p\}$,
\begin{equation}\label{eq:constrProb}
{\rm Pr} \left[\displaystyle \sum_{t=0}^{T-1}G_i(x^t,\xi_t) \leq \pi(T,\eta)\right] \geq 1-\eta,
\end{equation}
where
\begin{equation}\label{eq:thetadef}
\begin{array}{ll}
\pi (T,\eta)=&\kappa_3\left[1+2\sqrt{p}\kappa_g^2\right] +2\sqrt{p}\nu_g\kappa_g^2\\[6pt]
& \quad +[\left(1+2\sqrt{p}\kappa_g^2\right) (\kappa_0+\kappa_1+\kappa_4)+ 4\kappa_g\kappa_f]T^{1/2}\\[8pt]
&\quad +8\left(1+2\sqrt{p}\kappa_g^2\right)\displaystyle
\frac{\nu_g^2}{\epsilon_0}T^{1/2}\log \left( \displaystyle
\frac{T+1}{\eta}\right).
\end{array}
\end{equation}
\end{theorem}
{\bf Proof}. Define $Z(t)=\|\lambda^t\|$ for $\forall t \in
\{0,1,2,\ldots\}$. From Lemma \ref{Lemma 7}, $Z(t)$ satisfies the
conditions in Lemma \ref{Lemma 5} with $\delta_{\max}=\sigma \nu_g$,
$t_0=s$ and $\zeta=\displaystyle \frac{1}{2}\epsilon_0\sigma$ and
$$
 \vartheta =\displaystyle \frac{\epsilon_0\sigma s}{2}+ \nu_g\sigma (s-1)+\displaystyle \frac{\alpha D_0^2}{\epsilon_0s}+\displaystyle \frac{2\nu_f}{\epsilon_0}+
 \displaystyle \frac{\sigma \nu_g^2}{\epsilon_0}.
$$
Let $T\geq 1$ and $\lambda \in (0,1)$. Taking $\mu=\displaystyle
\frac{\eta}{T+1}$ in part 2 of  Lemma \ref{Lemma 5}, we obtain
\begin{equation}\label{eq:p1}
{\rm Pr}\,\left[\|\lambda^t\|\geq \gamma
(\sigma,\alpha,s,\eta)\right] \leq \displaystyle \frac{\eta}{T+1},\,
\forall t \in \{0,1,2,\ldots, T\},
\end{equation}
where
$$
\begin{array}{ll}
\gamma (\sigma, \alpha,s,\eta)& =\phi \left(\sigma, \alpha,s,\displaystyle \frac{\eta}{T+1}\right)\\[8pt]
&=\kappa_0+\kappa_1\displaystyle \frac{\alpha}{s}+\kappa_3
 \sigma+\kappa_4 \sigma s+8\displaystyle \frac{\nu_g^2}{\epsilon_0}\log \left( \displaystyle \frac{T+1}{\eta}\right)\sigma s.
 \end{array}
$$
This implies
$$
{\rm Pr}\,\left[\|\lambda^t\|\geq \gamma (\sigma, \alpha,s,\eta)
\mbox{ for some } t \in \{0,1,2,\ldots, T\}\right] \leq  \eta
 $$
or
\begin{equation}\label{eq:p2}
{\rm Pr}\,\left[\|\lambda^t\|\leq \gamma (\sigma,\alpha, s,\eta)
\mbox{ for } \forall t \in \{0,1,2,\ldots, T\} \right] \geq 1-\eta.
\end{equation}
It follows from (\ref{eq:c2}) and Assumption (A4)  that

\begin{equation}\label{eq:c2a}
\begin{array}{ll}
\displaystyle \sum_{t=0}^{T-1} G_i(x^t,\xi_t) &\leq \displaystyle
\frac{1}{\sigma} \lambda^{T}_i+  \displaystyle
\frac{2\kappa_g\kappa_f}{2\alpha-p\kappa_g^2\sigma }T +
\displaystyle \frac{\sqrt{p}\kappa_g^2}{2\alpha-p\kappa_g^2\sigma
}\displaystyle \sum_{t=0}^{T-1}\|\lambda^t\| +\displaystyle
\frac{\sqrt{p}\nu_g\kappa_g^2\sigma}{2\alpha-p\kappa_g^2\sigma }T.
\end{array}
\end{equation}
for $i=1,\ldots,p$. Thus, for $\sigma=T^{-1/2},\alpha=T^{1/2}$,
$2\sqrt{T}-p\kappa_g^2\sigma> \sqrt{T}/2$ when $T$ is very large, we
have from (\ref{eq:c2a}) that
\begin{equation}\label{eq:consV100a}
\begin{array}{ll}
\displaystyle \sum_{t=0}^{T-1} G_i(x^t,\xi_t) &\leq \displaystyle
T^{1/2}\|\lambda^{T}\|+
4\kappa_g\kappa_fT^{1/2}+2\sqrt{p}\nu_g\kappa_g^2
+2\sqrt{p}\kappa_g^2T^{-1/2} \sum_{t=0}^{T-1}\|\lambda^t\|.
\end{array}
\end{equation}
Let $s=\lceil T^{1/2} \rceil$. Noting, from (\ref{eq:consV100a}),  that
$$
\begin{array}{ll}
\displaystyle \sum_{t=0}^{T-1} G_i(x^t,\xi_t)
&\leq\left[1+2\sqrt{p}\kappa_g^2\right] T^{1/2} \gamma(T^{-1/2},T^{1/2},T^{1/2},\eta)+ 4\kappa_g\kappa_fT^{1/2}+2\sqrt{p}\nu_g\kappa_g^2\\[10pt]
&=\kappa_3\left[1+2\sqrt{p}\kappa_g^2\right] +2\sqrt{p}\nu_g\kappa_g^2+[\left(1+2\sqrt{p}\kappa_g^2\right) (\kappa_0+\kappa_1+\kappa_4)+ 4\kappa_g\kappa_f]T^{1/2}\\
&\quad +8\left(1+2\sqrt{p}\kappa_g^2\right)\displaystyle \frac{\nu_g^2}{\epsilon_0}T^{1/2}\log \left( \displaystyle \frac{T+1}{\eta}\right)\\
&=  \pi (T, \eta),
\end{array}
$$
when $\|\lambda^t\|\leq \gamma (T^{-1/2},T^{1/2},T^{1/2},\eta)$ for $\forall t \in \{0,1,2,\ldots, T\}$, we obtain the probability inequality (\ref{eq:constrProb}) from (\ref{eq:p2}). \hfill $\Box$\\
Define
$$
\omega_c (T)=\pi (T,e^{-T^{1/4}})/T,
$$
then
$$
\omega_c (T)=\left[8\left(1+2\sqrt{p}\kappa_g^2\right)\displaystyle
\frac{\nu_g^2}{\epsilon_0}\right] T^{-1/4}+{\rm o}(T^{-1/4}).
$$
We can obtain the following result from Theorem \ref{th:constr3} directly.
\begin{corollary}\label{th:constr3cor}  Let $(x^t,\lambda^t)$ be generated by PMMSopt, and Assumptions (A1)--(A4) be satisfied. If $\sigma=T^{-1/2},\, \alpha=T^{1/2}$ in PMMSopt, then for all $i\in \{1,\ldots,p\}$,
\begin{equation}\label{eq:constrProbcor}
{\rm Pr} \left[\displaystyle \frac{1}{T}\sum_{t=0}^{T-1}G_i(x^t,\xi_t) \leq \omega_c (T)\right] \geq 1-e^{-T^{1/4}}.
\end{equation}
\end{corollary}
For $\widehat x \in \Phi$, define $\widehat Z(0)=0$ and
$$\widehat Z(t)=\displaystyle \sum_{l=0}^{t-1} \langle \lambda^l, G(\widehat x, \xi_l) \rangle.$$
 Recall $\xi_{[0]}=\{\emptyset, \Omega\}$ and $\xi_{[t]}=\sigma (\xi_1,\ldots,\xi_t)$. In the following lemma, we will show that for any $c>0$, $\widehat Z(t)$ satisfies  conditions in Lemma 9 with ${\cal F}(t)=\xi_{[t]}$ and $Y(t)=\|\lambda^{t+1}\|-c/\nu_g$.

\begin{lemma}\label{lem:l13}
Let $\widehat x \in \Phi$. Let $c>0$  be arbitrary.  Let
$(x^t,\lambda^t)$ be generated by PMMSopt and Assumption (A3) be
satisfied. Define $\widehat Z(0)=0$ and $\widehat Z(t)=\displaystyle
\sum_{l=0}^{t-1} \langle \lambda^l, G(\widehat x, \xi_l) \rangle$,
$\forall t \geq 1$, then $\{\widehat Z(t),t \geq 0\}$ is a
supermartingale adapted to filtration $\{\xi_{[t]}, t\geq 0\}$ such
that
$$
\{|\widehat Z(t+1)-\widehat Z(t)|>c\} \subseteq \{ Y(t)> 0\}, \forall t \geq 0
$$
where $Y(t)=\|\lambda^{t+1}\|-c/\nu_g$ is a random variable adapted
to $\xi_{[t]}$.
\end{lemma}
{\bf Proof}. It is very easy to check $\{\widehat Z(t),t\geq 0\}$ is adapted to $\{\xi_{[t]}, t\geq 0\}$. Now we prove that
$\{\widehat Z(t),t \geq 0\}$ is a supermartingale. Since
$$
\widehat Z(t+1)=\widehat Z(t)+\langle \lambda^{t+1}, G(\widehat x, \xi_{t+1})\rangle,
$$
we have
$$
\begin{array}{ll}
\mathbb{E}[\widehat Z(t+1)\,|\, \xi_{[t]}]&=\mathbb{E}[\widehat Z(t)+\langle \lambda^{t+1}, G(\widehat x, \xi_{t+1})\rangle\,|\, \xi_{[t]}]\\[8pt]
& = \widehat Z(t)+\langle \lambda^{t+1}, \mathbb{E}[G(\widehat x, \xi_{t+1})]\rangle\\[8pt]
&=\widehat Z(t)+\langle \lambda^{t+1}, g(\widehat x)\rangle\\[8pt]
&\leq \widehat Z(t),
\end{array}
$$
which follows from $\widehat Z(t) \in \xi_{[t]}, \lambda^{t+1}\in \xi_{[t]}$, $G(\widehat x, \xi_{t+1})$ is independent of $\xi_{[t]}$ and $g(\widehat x)\leq0$. Thus we obtain that $\{\widehat Z(t),t\geq 0\}$ is a supermartingale.

From Assumption (A3), we get
$$
|\widehat Z(t+1)-\widehat Z(t)|=|\langle \lambda^{t+1}, G(\widehat x, \xi_{t+1} \rangle| \leq \nu_g\|\lambda^{t+1}\|.
$$
This implies that $\|\lambda^{t+1}\|> c/\nu_g$ if $|\widehat Z(t+1)-\widehat Z(t)| >c$ and
$$
\{|\widehat Z(t+1)-\widehat Z(t)| >c\} \subseteq \{\|\lambda^{t+1}\|> c/\nu_g\}.
$$
Since $\lambda^{t+1}$ is adapted to $\xi_{[t]}$, we have that $Y(t)=\|\lambda^{t+1}\|-c/\nu_g$ is a random variable adapted to
$\xi_{[t]}$. \hfill $\Box$\\

Next, we will use (\ref{eq:n13}) and Lemma \ref{Lemma 9} to establish a high probability objective regret bound.
\begin{theorem}\label{th:obj3} Let $\eta\in (0,1)$ and $\widehat x \in \Phi$.  Let $(x^t,\lambda^t)$ be generated by PMMSopt, and Assumption (A1)--Assumption (A4) be satisfied. If $\sigma=T^{-1/2},\, \alpha=T^{1/2}$ in PMMSopt, then
\begin{equation}\label{eq:ObjProb}
{\rm Pr} \left[\displaystyle \sum_{t=0}^{T-1}F(x^t,\xi_t) \leq
\displaystyle \sum_{t=0}^{T-1}F(\widehat x,\xi_t) + \beta
(T,\eta)\right] \geq 1-\eta,
\end{equation}
where
\begin{equation}\label{eq:betadef}
\begin{array}{ll}
\beta(T,\eta)=& \displaystyle \frac{\kappa_f^2+\nu_g^2}{2} T^{1/2}
+\displaystyle \frac{D_0^2}{2}T^{1/2}\\[10pt]
&\quad +\sqrt{2}\nu_g\log^{1/2}\left( \displaystyle \frac{2}{\eta}\right)\displaystyle\left[
(\kappa_0+\kappa_1+\kappa_4) T^{1/2}+\kappa_3\right.\\[10pt]
& \quad \left.+\displaystyle \frac{8\nu_g^2}{\epsilon_0}T^{1/2}\log\left( \displaystyle \frac{2T}{\eta}\right)
\right].
\end{array}
\end{equation}
\end{theorem}
{\bf Proof}. By Lemma \ref{lem:l13}, we know that $\widehat Z(t)$ satisfies conditions in Lemma \ref{Lemma 9}. Fix $T>0$, we obtain from Lemma \ref{Lemma 9} for any $c>0$ that
\begin{equation}\label{eq:e31}
{\rm Pr}\left[\displaystyle \sum_{t=0}^{T-1} \langle \lambda^t,
G(\widehat x, \xi_t)\rangle  \geq \gamma \right] \leq
e^{-\gamma^2/(2Tc^2)}+\displaystyle \sum_{t=0}^{T-1} {\rm Pr}\left[
\|\lambda^{t+1}\|> c/\nu_g\right].
\end{equation}
For given $\eta \in (0,1)$, we shall show how to choose $\gamma$ and $c$ such that each term in (\ref{eq:e31}) is not larger than $\eta/2$.

Noting that by Lemma \ref{Lemma 7}, random process
$Z(t)=\|\lambda^t\|$ satisfies conditions in Lemma \ref{Lemma 5}
with $\delta_{\max}= \nu_g \sigma$, $t_0=s$,  $\zeta=\displaystyle
\frac{\epsilon_0}{2} \sigma$ and
$$
 \vartheta =\displaystyle \frac{\epsilon_0\sigma s}{2}+ \nu_g
 \sigma (s-1)+\displaystyle \frac{\alpha D_0^2}{\epsilon_0s}+\displaystyle \frac{2\nu_f}{\epsilon_0}+
 \displaystyle \frac{\sigma \nu_g^2}{\epsilon_0}.
 $$
 For the second term being not larger than $\eta/2$, it suffices to choose $c$ such that
 $$
 {\rm Pr}\left[ \|\lambda^t\|>c/\nu_g\right] \leq \displaystyle \frac{\eta}{2T}, \,\forall t \in \{1,2,\ldots, T\}.
 $$
 The above inequality holds from  part 2 of  Lemma \ref{Lemma 5}  when we choose
 \begin{equation}\label{eq:cdef}
 c=c(\sigma,\alpha,s)=\left[\kappa_0+\kappa_1\displaystyle \frac{\alpha}{s}+\kappa_3
 \sigma+\kappa_4 \sigma s+\displaystyle \frac{8\nu_g^2}{\epsilon_0}\log \left( \displaystyle \frac{2T}{\eta}\right)\sigma s\right]\nu_g,
 \end{equation}
 where $s$ is an arbitrary integer. Define
 \begin{equation}\label{eq:gammc}
  \gamma (\sigma,\alpha,s,\eta) =\sqrt{2T} \log^{1/2} \left(\displaystyle \frac{2}{\eta}\right)c(\sigma,\alpha,s).
   \end{equation}
   Then, for $\gamma = \gamma (\sigma,\alpha,s,\eta)$ in (\ref{eq:e31}), the first term in this equation is equal to $\eta/2$. Thus we have, for $c = c(\sigma,\alpha,s)$ and $\gamma = \gamma(\sigma,\alpha,s,
   \eta)$ defined by (\ref{eq:cdef}) and (\ref{eq:gammc}), respectively, that
   $$
   {\rm Pr}\left[\displaystyle \sum_{t=0}^{T-1} \langle \lambda^t,G(\widehat x,\xi_t)\rangle\geq \gamma (\sigma,\alpha,s,\eta)\right]\leq \eta
   $$
   or equivalently
   \begin{equation}\label{eq:512}
   {\rm Pr}\left[\displaystyle \sum_{t=0}^{T-1} \langle \lambda^t,G(\widehat x,\xi_t)\rangle\leq \gamma (\sigma,\alpha,s,\eta)\right]\geq 1-\eta
         \end{equation}
    It follows from (\ref{eq:n13}) that
\begin{equation}\label{eq:514}
    \begin{array}{ll}
\displaystyle \sum_{t=0}^{T-1} F(x^t,\xi_t)  &
 \leq
\displaystyle \sum_{t=0}^{T-1} F(\widehat x,\xi_t)  + \displaystyle
\frac{1}{2\alpha}\kappa_f^2 T+\displaystyle \frac{\sigma}{2}
\nu_g^2T +\displaystyle \sum_{t=0}^{T-1} \langle
\lambda^t,G(\widehat x,\xi_t)\rangle
+\displaystyle \frac{\alpha}{2}D_0^2\\[12pt]
& = \displaystyle \sum_{t=0}^{T-1} F(\widehat x,\xi_t)  +
\displaystyle \frac{1}{2}(\kappa_f^2+\nu_g^2) T^{\frac{1}{2}}+
\displaystyle \sum_{t=0}^{T-1} \langle \lambda^t,G(\widehat
x,\xi_t)\rangle +\displaystyle \frac{D_0^2}{2}T^{\frac{1}{2}}.
\end{array}
\end{equation}
Taking $s=\lceil T^{1/2} \rceil$, we obtain from (\ref{eq:514}) that if
$$
\displaystyle \sum_{t=0}^{T-1} \langle \lambda^t,G(\widehat
x,\xi_t)\rangle\leq \gamma (T^{-1/2},T^{1/2},T^{1/2},\eta),
$$
then
$$
 \begin{array}{ll}
\displaystyle \sum_{t=0}^{T-1} F(x^t,\xi_t)  &
 \leq
\displaystyle \sum_{t=0}^{T-1} F(\widehat x,\xi_t)  + \displaystyle
\frac{1}{2}(\kappa_f^2+\nu_g^2) T^{1/2} \\[10pt]
& \quad +\gamma (T^{-1/2},T^{1/2},T^{1/2},\eta)
+\displaystyle \frac{D_0^2}{2}T^{1/2}\\[10pt]
& \leq \displaystyle \sum_{t=0}^{T-1} F(\widehat x,\xi_t)+
\displaystyle \frac{1}{2}(\kappa_f^2+\nu_g^2) T^{1/2}
+\displaystyle \frac{D_0^2}{2}T^{1/2}\\[10pt]
&\quad +\sqrt{2}\nu_g\log^{1/2}\left( \displaystyle \frac{2}{\eta}\right)\displaystyle\left[
(\kappa_0+\kappa_1+\kappa_4) T^{1/2}+\kappa_3\right.\\[10pt]
& \quad \left.+\displaystyle \frac{8\nu_g^2}{\epsilon_0}T^{1/2}\log\left( \displaystyle \frac{2T}{\eta}\right)
\right]\\[10pt]
&= \displaystyle \sum_{t=0}^{T-1} F(\widehat x,\xi_t)
+\beta(T,\eta).
\end{array}
$$
We obtain the probability bound (\ref{eq:ObjProb}) from (\ref{eq:512}). \hfill $\Box$\\
Define
\begin{equation}\label{eq:wo}
\omega_o (T)=\beta (T,e^{-T^{1/4}})/T,
\end{equation}
then
$$
\omega_o (T)=\displaystyle \frac{8\sqrt{2}\nu_g^3}{\epsilon_0} T^{-1/8}+{\rm o}(T^{-1/8}).
$$
We can obtain the following result from Theorem \ref{th:obj3} directly.
\begin{corollary}\label{th:obj3cor} Let $\eta\in (0,1)$ and $\widehat x \in \Phi$.  Let $(x^t,\lambda^t)$ be generated by PMMSopt, and Assumptions (A1)--(A4) be satisfied. If $\sigma=T^{-1/2},\, \alpha=T^{1/2}$ in PMMSopt, then
\begin{equation}\label{eq:ObjProbcor}
{\rm Pr} \left[\displaystyle \displaystyle
\frac{1}{T}\sum_{t=0}^{T-1}F(x^t,\xi_t) \leq
\displaystyle\displaystyle \frac{1}{T} \sum_{t=0}^{T-1}F(\widehat
x,\xi_t) + \omega_o (T)\right] \geq 1-e^{-T^{1/4}}.
\end{equation}
\end{corollary}

\section{Conclusion}\label{Sec5}
\setcounter{equation}{0}
In this paper, for the first time we present a stochastic approximation proximal method of multipliers (PMMSopt) for solving convex stochastic programming with expectation constraints. We  show that, when the objective and constraint functions are generally convex,  this algorithm exhibits  ${\rm O}(T^{-1/2})$  objective regret and   ${\rm O}(T^{-1/2})$ constraint violation regret if parameters in the algorithm are properly chosen, where $T$ denotes the number of iterations.  Moreover, we show that, with at least $1-e^{-T^{1/4}}$ probability, the algorithm has no more than  ${\rm O}(T^{-1/4})$ objective regret and no more than  ${\rm O}(T^{-1/8})$ constraint violation  regret.

The research
presented in this paper has provided a proximal point-type method for solving convex  stochastic programming with expectation constraints, which  has the complexity bound  comparable to the currently best existing one in  \cite{LanZ2016,YMNeely2017}.  However, there are several questions to answer for  the proposed stochastic approximation proximal method of multipliers.
 One question that has not been answered here is how to obtain the same complexity bound when $X_0$ in the problem is unbounded.
  The complexity bound in this paper requires that $x^{t+1}$ in Step 1 of PMMSopt is an  exact solution  to the subproblem, so the  the second question is whether the same complexity bound can be derived when $x^{t+1}$ is inexactly solved.   The algorithms in \cite{LanZ2016,YMNeely2017} only require solving simple optimization
problems involving the stochastic gradients of $f$ and $g_i$'s, the third question is whether we can construct a proximal point-type method in which subproblems for $x^{t+1}$ are simple optimization problems, or even have explicit solutions.
\mbox{}\\[10pt]

\end{document}